\documentclass[a4paper,fleqn]{cas-dc}

\usepackage[numbers]{natbib}
\usepackage{graphicx}
\usepackage{amsmath}
\usepackage{bm}
\usepackage{amssymb}
\usepackage{algorithm,algorithmic}
\usepackage{graphicx}
\usepackage{fancyhdr}
\usepackage{amsmath}
\usepackage{mathtools}
\usepackage{setspace}
\usepackage{algorithmic}
\usepackage{algorithm}
\usepackage{nameref}
\usepackage{optidef}
\usepackage{amssymb}
\usepackage{amsfonts}
\usepackage{bm}
\usepackage{textcomp}

\usepackage{nameref}
\usepackage{subcaption}
\usepackage{url}

\def\tsc#1{\csdef{#1}{\textsc{\lowercase{#1}}\xspace}}
\tsc{WGM}
\tsc{QE}
\tsc{EP}
\tsc{PMS}
\tsc{BEC}
\tsc{DE}

\begin{document}

\let\WriteBookmarks\relax
\def\floatpagepagefraction{1}
\def\textpagefraction{.001}

\shorttitle{Container pre-marshalling problem minimizing CV@R under uncertainty of ship arrival times}

\shortauthors{Ikuma et~al.}

\title [mode = title]{Container pre-marshalling problem minimizing CV@R under uncertainty of ship arrival times}                      



%
\author[1]{Daiki Ikuma}





\affiliation[1]{organization={Graduate School of Science and Technology, University of Tsukuba},
    addressline={1--1--1 Tennodai}, 
    city={Tsukuba-shi},
    postcode={305-8573}, 
    state={Ibaraki},
    country={Japan}}

\author[1]{Shunnosuke Ikeda}
\cormark[1]

\author[2]{Noriyoshi Sukegawa}


\affiliation[2]{organization={Department of Advanced Sciences, Faculty of Science and Engineering, Hosei University},
    addressline={3--7--2 Kajinocho}, 
    city={Koganei-shi},
    postcode={184--8584}, 
    state={Tokyo},
    country={Japan}}

\author[3]{Yuichi Takano}
\cormark[2]

\affiliation[3]{organization={Institute of Systems and Information Engineering, University of Tsukuba},
    addressline={1--1--1 Tennodai}, 
    city={Tsukuba-shi},
    postcode={305--8573}, 
    state={Ibaraki},
    country={Japan}}

\cortext[cor1]{Corresponding author}
\cortext[cor2]{Principal corresponding author}



\begin{abstract}
This paper is concerned with the container pre-marshalling problem, which involves relocating containers in the storage area so that they can be efficiently loaded onto ships without reshuffles. 
In reality, however, ship arrival times are affected by various external factors, which can cause the order of container retrieval to be different from the initial plan. 
To represent such uncertainty, we generate multiple scenarios from a multivariate probability distribution of ship arrival times. 
We derive a mixed-integer linear optimization model to find an optimal container layout such that the conditional value-at-risk is minimized for the number of misplaced containers responsible for reshuffles. 
Moreover, we devise an exact algorithm based on the cutting-plane method to handle large-scale problems. 
Numerical experiments using synthetic datasets demonstrate that our method can produce high-quality container layouts compared with the conventional robust optimization model. 
Additionally, our algorithm can speed up the computation of solving large-scale problems. 
\end{abstract}



\begin{keywords}
container terminal \sep  pre-marshalling \sep conditional value-at-risk (CV@R) \sep  cutting-plane method \sep mixed-integer optimization 
\end{keywords}

\maketitle

\section{INTRODUCTION}
\label{sec:intro}
\subsection{Background}
Maritime transportation is the most cost-effective means of transporting large quantities of goods and raw materials around the world, accounting for over 90\% of the total volume of world trade~\cite{ICS2019}.
In particular, containerized transportation accounts for about 60\% of the total volume of global maritime trade and thus plays a significant role as a means of transportation~\cite{ICS2020}.
In addition, the volume of containerized transportation has been increasing annually~\cite{MarineInsight}. 
This leads to a growing demand for more efficient maritime transportation, especially in terms of operational efficiency at container terminals that connect land and maritime transportation.

Containers are usually stacked for storage due to limited space at container terminals. 
Specifically, containers are placed on the ground first, and once the space on the ground is filled, additional containers are stacked on top of other containers.
This stacking allows us to make the most of the available space and increase storage capacity, but it also increases the complexity of retrieving a particular container from a stack of stored containers.

If a container stacked below needs to be retrieved, all containers above it must be moved to other stacks. 
Such unproductive movement of containers is referred to as the \emph{reshuffle}.
Since cranes or other machines are required to move containers, the number of reshuffles is a critical factor in increasing the efficiency of loading containers onto ships.
Unnecessary movement of containers also causes various problems such as maritime congestion and increased fuel costs for cranes~\cite{maguire2010relieving,goodchild2005crane}.

To improve loading efficiency at container terminals, it is effective to relocate containers in advance so that they will be retrieved without reshuffles. 
In reality, however, a ship arrival time is affected by a variety of external factors, such as bad weather, current variation, cargo handling delays, and ship malfunction. 
These factors can cause the order of container retrieval to be different from the initial plan, so we should take into account the uncertainty of ship arrival times to develop an effective container pre-marshalling strategy.

\subsection{Related work}
The problem of relocating containers so that they can be retrieved without reshuffles is known as the \emph{container pre-marshalling problem}~\cite{lee2007optimization}. 
A number of studies have been conducted on container pre-marshalling problems; see also comprehensive surveys~\cite{lehnfeld2014loading,caserta2020container,lersteau2022survey}. 
The container pre-marshalling problem can be formulated as an optimization model, and various exact and heuristic algorithms have been proposed~\cite{bortfeldt2012tree,huang2012heuristic,tierney2017solving,tanaka2018solving,parreno2019integer}. 
Methods using constraint programming~\cite{jimenez2023constraint} and deep learning~\cite{hottung2020deep} have also gained attention in recent years.
These studies assume that all ships arrive as scheduled.

Several studies have considered the uncertainty of ship arrival times in the container pre-marshalling problem~\cite{rendl2013constraint,tierney2016solving,boywitz2018robust,boge2020robust,zweers2020optimizing,maniezzo2021stochastic}.
Rendl et al.~\cite{rendl2013constraint} proposed modelling arrival times as time intervals and relocating containers so that these intervals do not overlap in every stack. 
Tierney et al.~\cite{tierney2016solving} generalized this model by introducing a blocking matrix that specifies which containers can be stacked on top of others, and provided an efficient method for solving this generalized problem. 
However, these problems can be infeasible, meaning that these models could fail to output any suggestion for how to relocate containers. 

Boywitz and Boysen~\cite{boywitz2018robust} proposed to maximize the minimum margin between the ship arrival times of any two consecutive containers stored in the same stack. 
Boge et al.~\cite{boge2020robust} formulated a robust optimization model to minimize the number of misplaced containers responsible for reshuffles based on an uncertainty set of ship arrival orders. 
However, it would be desirable to incorporate more detailed information about the uncertainty of ship arrival times as a multivariate probability distribution. 

Zweers et al.~\cite{zweers2020optimizing} proposed to minimize the expected number of reshuffles for a specific policy. 
Maniezzo et al.~\cite{maniezzo2021stochastic} proposed to minimize the expected cost of reshuffles, where the uncertainty of ship arrival times is modeled using a Poisson distribution. 
These models focus only on expected values and do not control the risk of suffering a large number of reshuffles.

\subsection{Our contribution}
To overcome these challenges posed in the prior studies, we adopt the conditional value-at-risk (CV@R)~\cite{rockafellar2002conditional}, a widely used risk measure in various application domains~\cite{filippi2020conditional}. 
We generate multiple scenarios from a multivariate probability distribution to represent the uncertainty of ship arrival times. 
We then derive a mixed-integer linear optimization model to find an optimal container layout such that the conditional value-at-risk is minimized for the number of misplaced containers responsible for reshuffles. 
However, since the size of this optimization problem is highly dependent on the number of scenarios, it is very difficult to solve the problem in a practical time frame when a large number of scenarios are considered. 
To handle such large-scale problems, we devise an exact algorithm based on the cutting-plane method~\cite{kunzi2006computational,takano2015cutting,kobayashi2021bilevel}. 

Numerical experiments using synthetic datasets demonstrate that our method can produce high-quality container layouts under uncertainty of ship arrival times compared with the robust optimization model proposed in Boge et al.~\cite{boge2020robust}. 
Moreover, our algorithm can speed up the computation of solving the container pre-marshalling problem especially when the problem size is large. 

\subsection{Structure of this paper}
This paper is structured as follows:
Section~\ref{sec:prob} describes the container pre-marshalling problem under uncertainty of ship arrival times.
Section~\ref{sec:form} formulates our mixed-integer linear optimization model to minimize CV@R for the number of misplaced containers. 
Section~\ref{sec:algo} details our solution algorithm based on the cutting-plane method. 
Section~\ref{sec:exper} evaluates the effectiveness of our method based on numerical results. 
Section~\ref{sec:concl} concludes with a brief summary of this study and a discussion of future research directions.

\section{PROBLEM DESCRIPTION}
\label{sec:prob}
In this section, we first describe the container pre-marshalling problem and then discuss the motivation for considering uncertainty of ship arrival times.

\subsection{Notation}
\begin{figure*}[t]
\centering
    \begin{minipage}[t]{0.24\linewidth}
        \centering
        \includegraphics[keepaspectratio, width=0.90\linewidth]{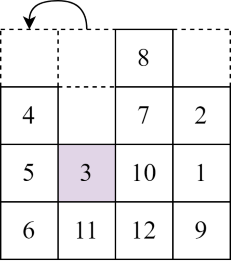}
        \subcaption{}\label{fig:pm_a}
    \end{minipage}%
    \begin{minipage}[t]{0.24\linewidth}
        \centering
        \includegraphics[keepaspectratio, width=0.90\linewidth]{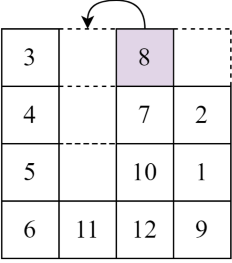}
        \subcaption{}\label{fig:pm_b}
    \end{minipage}%
    \begin{minipage}[t]{0.24\linewidth}
        \centering
        \includegraphics[keepaspectratio, width=0.90\linewidth]{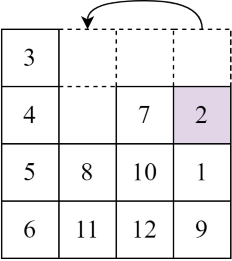}
        \subcaption{}\label{fig:pm_c}
    \end{minipage}%
    \begin{minipage}[t]{0.24\linewidth}
        \centering
        \includegraphics[keepaspectratio, width=0.90\linewidth]{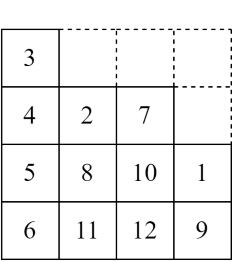}
        \subcaption{}\label{fig:pm_d}
    \end{minipage}
\caption{Example of moving containers in pre-marshalling}
\label{fig:figure_label}
\end{figure*}

A container storage area, also known as a bay, consists of $S$ stacks arranged in a single line, and the height of each stack is limited to $H$ containers. 
All containers are assumed to have the same size, so the storage area is two-dimensional, as shown in Figure~\ref{fig:figure_label} ($S=4$ and $H=4$).

Let $\mathcal{S} \coloneqq \{1,2,\ldots,S\}$ and $\mathcal{H} \coloneqq \{1,2,\ldots,H\}$ represent the index sets of stacks (columns) and height levels (rows), respectively.
Any position in the storage area can be denoted as $(s, h) \in \mathcal{S} \times \mathcal{H}$.
For example, the bottom left position is $(s,h) = (1,1)$, and the top right position is $(s,h) = (4,4)$ in Figure~\ref{fig:figure_label}. 

We also define the index set of containers as $\mathcal{N} \coloneqq \{1, 2, \ldots, N\}$ and that of container priority classes as $\mathcal{R} \coloneqq \{1, 2, \ldots, R\}$. 
The index set of containers in priority class $r \in \mathcal{R}$ is defined as $\mathcal{N}_r \subseteq \mathcal{N}$. 
For example, if $r < r'$, container $n \in \mathcal{N}_r$ is retrieved before container $n' \in \mathcal{N}_{r'}$. 
This priority class corresponds to the ship arrival order; that is, container $n \in \mathcal{N}_r$ is loaded onto the ship that arrives $r$th at the terminal. 

\vspace{3mm}

\subsection{Container pre-marshalling problem}
If a container stacked below need to be retrieved, all containers above it must be \emph{reshuffled} (i.e., moved to other stacks). 
The \emph{container pre-marshalling problem} aims to relocate containers so that they can be retrieved without reshuffles. 

A container in priority class $r\in \mathcal{R}$ at position $(s,h) \in \mathcal{S} \times \mathcal{H}$ is \emph{misplaced} if there exists a container in priority class $r'$ at position $(s',h')$ that is placed below in the same stack but retrieved earlier (i.e., $s=s'$, $h > h'$, and $r > r'$). 
Note that misplaced containers inevitably cause reshuffles. 

Figure \ref{fig:figure_label} shows an example of moving containers in pre-marshalling, where the number written on each container is its priority class. 
The initial container layout (Figure \ref{fig:pm_a}) has misplaced containers at positions $(3, 4)$ and $(4, 3)$. 
A sequence of three container moves can transform them into an ideal layout (Figure \ref{fig:pm_d}) without any misplaced containers.

\begin{figure}[t]
    \begin{minipage}[t]{0.49\linewidth}
        \centering
        \includegraphics[keepaspectratio, scale=0.43]{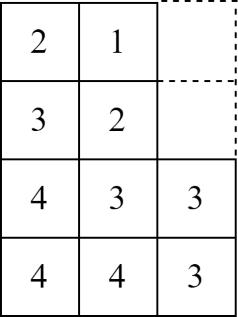}
        \subcaption{}\label{fig:cp_a}
    \end{minipage}
    \begin{minipage}[t]{0.49\linewidth}
      \centering
      \includegraphics[keepaspectratio, scale=0.43]{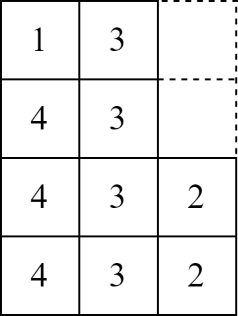}
      \subcaption{}\label{fig:cp_b}
    \end{minipage}
   \caption{Two exmaples of container layouts}
   \label{fig:figure_label2}
\end{figure}
\renewcommand\thesubfigure{\alph{subfigure}}

\subsection{Uncertainty of ship arrival times}
We introduce the uncertainty of ship arrival times into the container pre-marshalling problem. 
Following Boge et al.~\cite{boge2020robust}, we consider two layouts, neither of which contains misplaced containers, as shown in Figure~\ref{fig:figure_label2}. 

In Figure~\ref{fig:cp_a}, however, if the arrival order of the first and second ships is reversed, the container at position $(2, 4)$ will be misplaced.
Additionally, if the arrival order of the second and third ships is reversed, the containers at positions $(1, 4)$ and $(2, 3)$ will also be misplaced. 

In Figure~\ref{fig:cp_b}, on the other hand, a misplaced container will appear at position $(1, 4)$ only when the arrival order of the first and fourth ships is reversed.
This reversal is unlikely because these ships have the furthest scheduled arrival times.
For this reason, Figure~\ref{fig:cp_b} is more ideal than Figure~\ref{fig:cp_a} when ship arrival times are uncertain. 

\section{PROBLEM FORMULATION}
\label{sec:form}
In this section, we first give an overview of scenario generation and CV@R. 
We next formulate the container pre-marshalling problem minimizing CV@R for the number of misplaced containers. 

\subsection{Scenario generation}
\label{sec:sce}
To generate a set of scenarios for ship arrival orders, we sample instances of ship arrival times from a multivariate probability distribution.

Let $\mathcal{I}$ be the index set of sampled instances of ship arrival times. 
For each instance $i \in \mathcal{I}$, let 
\[
\bm{t}(i) \coloneqq (t_r(i))_{r \in \mathcal{R}}
\]
be a vector of ship arrival times, where $t_r(i)$ is the arrival time of the $r$th ship. 
An example of four instances with three ships is given as follows: 
\begin{align}
    & \bm{t}(1) = (0.3, 2.2, 2.5), \quad \bm{t}(2) = (1.9, 1.3, 2.6), \notag \\
    & \bm{t}(3) = (1.0, 1.8, 2.9), \quad \bm{t}(4) = (0.5, 2.8, 2.7). \notag
\end{align} 

Let $\mathcal{J}$ be the index set of scenarios for ship arrival orders. 
For each scenario $j \in \mathcal{J}$, let $p_j$ be the occurrence probability, and 
\[
\bm{o}(j) \coloneqq (o_r(j))_{r \in \mathcal{R}}
\]
be a vector representing the order of ship arrivals, where $o_r(j)$ indicates the ship that arrives $r$th at the terminal. 

From the example above, we have 
\begin{align}
    & \bm{o}(1) = (1, 2, 3), \quad p_1 = 1/2, \quad \because \bm{t}(1),\bm{t}(3) \notag \\
    & \bm{o}(2) = (2, 1, 3), \quad p_2 = 1/4, \quad \because \bm{t}(2) \notag \\
    & \bm{o}(3) = (1, 3, 2), \quad p_3 = 1/4. \quad \because \bm{t}(4) \notag
\end{align}
Here, $\bm{o}(1)$ means that all three ships arrive at the terminal as scheduled in the first scenario, whereas $\bm{o}(2)$ means that the second ship arrives at the terminal first in the second scenario. 

\begin{figure}[t]
\centering
\includegraphics[keepaspectratio,scale=1.2,width=\linewidth]
{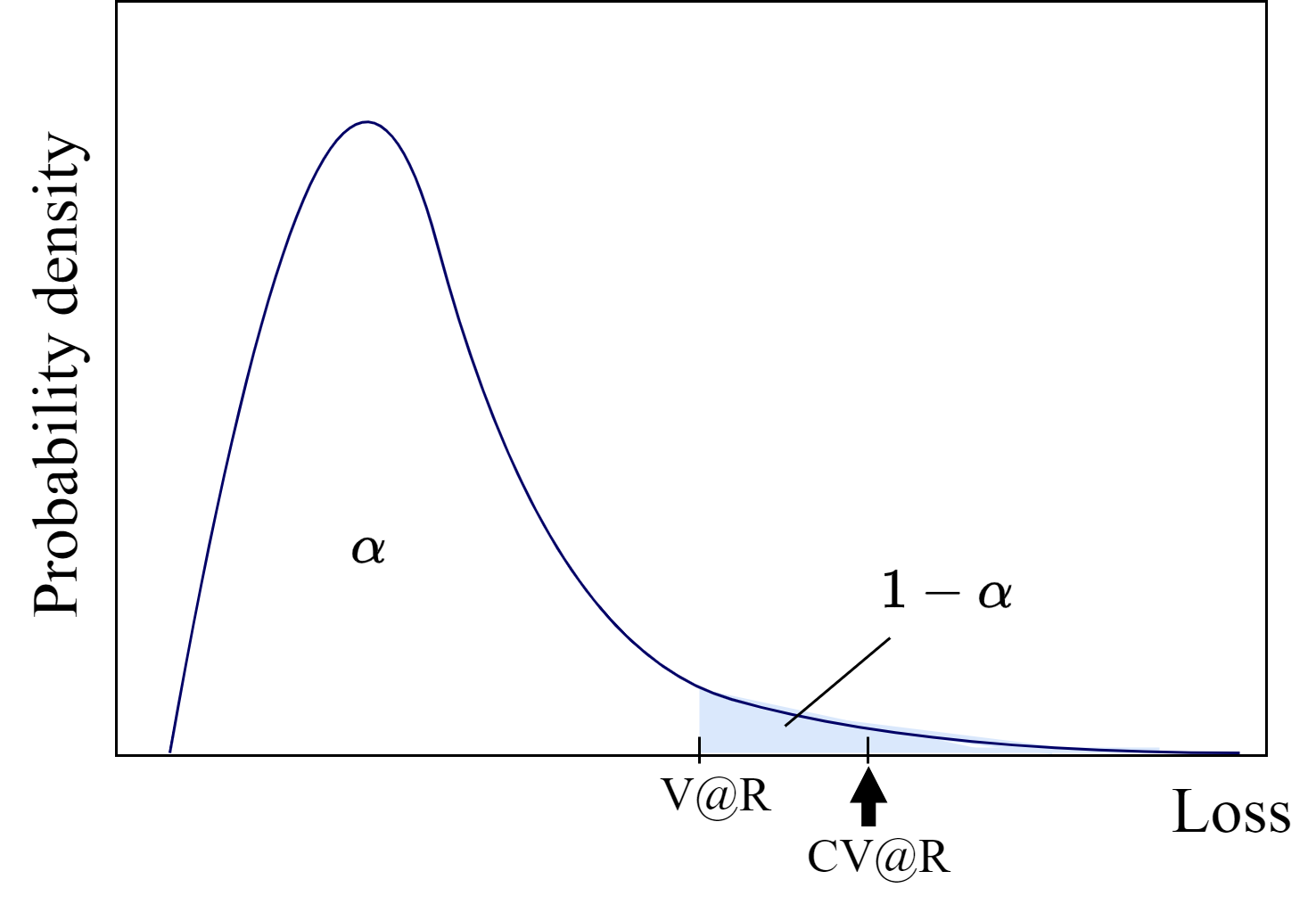}
\caption{Illustration of V@R and CV@R}
\label{fig:figure_label3}
\end{figure}

\subsection{Conditional value-at-risk (CV@R)}\label{subseq:CVaR}
As shown in Figure \ref{fig:figure_label3}, the value-at-risk (V@R) is the $\alpha$-quantile of the probability distribution of a random loss, where $\alpha \in [0, 1)$ is a user-defined probability level (e.g., $\alpha \in \{0.75, 0.90\}$). 
The conditional value-at-risk (CV@R) can be regarded as the conditional expectation of a random loss exceeding V@R.

Let $f(\bm{o}(j))$ be a loss function (e.g., the number of misplaced containers), which depends on the ship arrival order $\bm{o}(j)$ in scenario $j \in \mathcal{J}$. 
When V@R is denoted by $\gamma$, CV@R is written as follows~\cite{rockafellar2002conditional}: 
\[
\mbox{CV@R} \coloneqq \gamma + \frac{1}{1-\alpha} \sum_{j \in \mathcal{J}} p_j [f(\bm{o}(j)) - \gamma]_{+},
\]
where $[\,\cdot\,]_{+} \coloneqq \max\{\,\cdot\,, 0\}$. 

Let $\bm{v} \coloneqq (v_j)_{j \in J}$ be an auxiliary decision variable for representing the operation $[\,\cdot\,]_{+}$. 
CV@R can then be calculated using the lifting representation~\cite{rockafellar2002conditional,fabian2008handling}: 
\begin{align}
& \underset{{\bm v}, \gamma}{\text{min}} \quad \gamma + \frac{1}{(1-\alpha)}\sum_{j \in \mathcal{J}} p_j v_j \qquad \text{s.~t.} \label{obj:CVaR1} \\
& v_j\geq f(\bm{o}(j)) - \gamma \quad (j \in \mathcal{J}), \label{con1:CVaR1} \\
& v_j\geq 0 \quad (j \in \mathcal{J}), \label{con2:CVaR1}
\end{align}
where $\gamma$ is a decision variable corresponding to V@R. 

It is also known that CV@R can be calculated using the cutting-plane representation~\cite{kunzi2006computational,fabian2008handling}: 
\begin{align}
& \underset{u, \gamma}{\text{min}} \quad \gamma + u \qquad \text{s.~t.} \label{obj:CVaR2} \\
& u \ge \frac{1}{(1-\alpha)}\sum_{j \in \mathcal{G}} p_j \left( f(\bm{o}(j)) -\gamma \right) \quad (\mathcal{G} \in 2^{\mathcal{J}}), \label{con1:CVaR2}
\end{align}
where $u$ is an auxiliary decision variable for calculating CV@R. 
Note that $\mathcal{G} \subseteq \mathcal{J}$ is a subset of scenarios drawn from the power set of $\mathcal{J}$.

\subsection{Decision variables}
We define decision variables to be used in the container pre-marshalling problem. 
Let 
\[
\bm{x} \coloneqq (x_{shr})_{(s,h,r) \in \mathcal{S} \times \mathcal{H} \times \mathcal{R}}
\]
be a binary decision variable for determining container positions, where $x_{shr} = 1$ if a container in priority class $r \in \mathcal{R}$ is placed at position $(s,h) \in \mathcal{S} \times \mathcal{H}$, and $x_{shr} = 0$ otherwise. 

Let 
\[
\bm{z} \coloneqq (z_{shj})_{(s,h,j) \in \mathcal{S} \times (\mathcal{H} \setminus \{1\}) \times \mathcal{J}}
\]
be a real-valued decision variable for identifying misplaced containers. In scenario $j \in \mathcal{J}$, $z_{shj} = 1$ if a container at position $(s,h) \in \mathcal{S} \times (\mathcal{H} \setminus \{1\})$ is misplaced, and $z_{shj} = 0$ otherwise. 

\subsection{Constraints}
\label{sec:const}
We formulate various constraints to be imposed on our optimization problem. 
The storage area must have a specified number of containers for each priority class $r \in \mathcal{R}$:
\begin{align}
\sum_{s\in \mathcal{S}}\sum_{h\in \mathcal{H}}x_{shr} = |\mathcal{N}_r| \quad (r \in \mathcal{R}). \label{eq:sum_x_Nr}
\end{align}

There should not be more than one container at any position $(s,h) \in \mathcal{S} \times \mathcal{H}$ within the storage area: 
\begin{align}
\sum_{r \in \mathcal{R}}x_{shr} \leq 1 \quad (s \in \mathcal{S}, ~h \in \mathcal{H}).\label{eq:sum_x_1}
\end{align}

If there is a container at position $(s,h+1) \in \mathcal{S} \times \mathcal{H}$, there must also be a container at the position immediately below it in the same stack: 
\begin{align}
& \sum_{r \in \mathcal{R}}x_{s,h+1,r} \leq \sum_{r \in \mathcal{R}}x_{shr} \notag \\
& \qquad (s \in \mathcal{S},~h \in \mathcal{H}\setminus\{H\}). \label{eq:sum_x_sum_x}
\end{align}

To count the number of misplaced containers, let us consider a situation where 
\[
\sum_{k=r}^R x_{shk} = 1, \quad \sum_{k=r}^R x_{sh'k} = 0
\]
for $h > h'$.
The first equation implies that there is a container at the position $(s,h)$, and its priority class is $r$ or more. 
The second equation implies that the priority class of the position $(s,h')$ is less than $r$.
Accordingly, the container at the position $(s,h)$ is misplaced.  

In light of this observation, misplaced containers for the ship arrival order $\bm{o}(j)$ can be identified as follows: 
\begin{align}
   & z_{shj}\geq \sum_{k=r}^R x_{sho_k(j)}-\sum_{k=r}^R x_{sh'o_k(j)} \notag\\
   & \qquad (s \in \mathcal{S},~(h,h') \in \mathcal{H}^*,~r \in \mathcal{R}, ~j \in \mathcal{J}), \label{eq:z_shj} \\
   & 0 \le z_{shj}\le 1 \quad (s\in \mathcal{S}, ~h\in \mathcal{H}\setminus\{1\}, ~j\in \mathcal{J}), \label{eq:z_shj_0}
\end{align}
where $\mathcal{H}^* \coloneqq \{(h,h') \in \mathcal{H} \times \mathcal{H} \mid h > h'\}$. 
Hence, the number of misplaced containers in each scenario $j \in \mathcal{J}$ can be calculated as follows: 
\[
f(\bm{o}(j)) = \sum_{s \in \mathcal{S}} \sum_{h \in \mathcal{H} \setminus \{1\}} z_{shj}. 
\]

\subsection{Mixed-integer linear optimization formulation} \label{subsubsec:formulation}

We formulate a container pre-marshalling problem to minimize CV@R for the number of misplaced containers.
We use the lifting representation of CV@R (Eqs.~\eqref{obj:CVaR1}--\eqref{con2:CVaR1}) with the constraints (Eqs.~\eqref{eq:sum_x_Nr}--\eqref{eq:z_shj_0}) given in Section~\ref{sec:const}. 
Our container pre-marshalling problem can be formulated as the following mixed-integer linear optimization problem: 
\begin{align}
    & \underset{\bm{v}, \bm{x}, \bm{z}, \gamma}{\text{min}} \quad \gamma + \frac{1}{1-\alpha}\sum_{j\in \mathcal{J}} p_j v_j \qquad \text{s.~t.} \label{obj:prob1}\\
    & v_j \geq \sum_{s\in \mathcal{S}}\sum_{h\in \mathcal{H}\setminus\{1\}} z_{shj}-\gamma \quad (j \in \mathcal{J}), \label{con1:prob1} \\
    & v_j \geq 0 \quad (j\in \mathcal{J}), \label{con2:prob1} \\
    & \sum_{s\in \mathcal{S}}\sum_{h\in \mathcal{H}}x_{shr} = |\mathcal{N}_r| \quad (r \in \mathcal{R}), \label{con3:prob1} \\
    & \sum_{r \in \mathcal{R}}x_{shr} \leq 1 \quad (s \in \mathcal{S}, ~h \in \mathcal{H}), \label{con4:prob1} \\
    & \sum_{r \in \mathcal{R}}x_{s,h+1,r} \leq \sum_{r \in \mathcal{R}}x_{shr} \notag \\
    & \qquad (s \in \mathcal{S}, ~h \in \mathcal{H}\setminus\{H\}),  \label{con5:prob1}\\
    & x_{shr} \in \{0, 1\} \quad (s \in \mathcal{S}, ~h \in \mathcal{H}, ~r \in \mathcal{R}), \label{con6:prob1} \\
    & z_{shj}\geq \sum_{k=r}^R x_{sho_k(j)}-\sum_{k=r}^R x_{sh'o_k(j)} \notag\\
    & \qquad (s \in \mathcal{S}, ~(h,h') \in \mathcal{H}^*,~r \in \mathcal{R}, ~j \in \mathcal{J}), \label{con7:prob1} \\
    & 0 \le z_{shj}\le 1 \quad (s\in \mathcal{S}, ~h\in \mathcal{H}\setminus\{1\}, ~j\in \mathcal{J}). \label{con8:prob1}
\end{align}

This problem can be solved exactly using a mathematical optimization solver.
By solving this problem, we can obtain an ideal container layout (e.g., Figure~\ref{fig:cp_b}) under uncertainty of ship arrival times. 

\section{SOLUTION ALGORITHM}
\label{sec:algo}
The size of problem~\eqref{obj:prob1}--\eqref{con8:prob1} is highly dependent on the number of scenarios, so it becomes difficult to solve this problem within a practical time frame as the number of scenarios increases.
To overcome this challenge, we develop an exact algorithm based on the cutting-plane method~\cite{kunzi2006computational,takano2015cutting,kobayashi2021bilevel}.

We focus on the following problem formulation based on the cutting-plane representation of CV@R~(Eqs.~\eqref{obj:CVaR2}--\eqref{con1:CVaR2}): 
\begin{align}
    & \underset{u, \bm{x}, \bm{z}, \gamma}{\text{min}} \quad \gamma + u \qquad \text{s.~t.} \label{obj:prob2}\\
    & u \ge \frac{1}{1-\alpha}\sum_{j \in \mathcal{G}} p_j \left( \sum_{s\in \mathcal{S}}\sum_{h\in \mathcal{H}\setminus\{1\}} z_{shj}-\gamma \right) \notag \\
    & \qquad (\mathcal{G} \in 2^{\mathcal{J}}), \label{con1:prob2}\\
    & z_{shj}\geq \sum_{k=r}^R x_{sho_k(j)}-\sum_{k=r}^R x_{sh'o_k(j)} \notag\\
    & \qquad (s \in \mathcal{S}, ~(h,h') \in \mathcal{H}^*,~r \in \mathcal{R}, ~j \in \mathcal{J}), \label{con2:prob2}\\
    & 0 \le z_{shj} \le 1 \quad (s\in \mathcal{S}, ~h\in \mathcal{H}\setminus\{1\}, ~j\in \mathcal{J}), \label{con3:prob2} \\
    & \mbox{Eqs.~\eqref{con3:prob1}--\eqref{con6:prob1}}. \label{con4:prob2}
\end{align}

We relax Eqs.~\eqref{con1:prob2}--\eqref{con3:prob2} using only a scenario subset $\mathcal{G}_0$ to obtain the relaxed problem:  
\begin{align}
    & \underset{u, \bm{x}, \bm{z}, \gamma}{\text{min}} \quad \gamma + u \qquad \text{s.~t.} \label{obj:prob3}\\
    & u \ge \frac{1}{1-\alpha}\sum_{j \in \mathcal{G}_0} p_j \left( \sum_{s\in \mathcal{S}}\sum_{h\in \mathcal{H}\setminus\{1\}} z_{shj}-\gamma \right), \label{con1:prob3}\\
    & z_{shj}\geq \sum_{k=r}^R x_{sho_k(j)}-\sum_{k=r}^R x_{sh'o_k(j)} \notag\\
    & \qquad (s \in \mathcal{S}, ~(h,h') \in \mathcal{H}^*,~r \in \mathcal{R}, ~j \in \mathcal{G}_0), \label{con2:prob3} \\
    & 0 \le z_{shj} \le 1 \quad (s\in \mathcal{S}, ~h\in \mathcal{H}\setminus\{1\}, ~j\in \mathcal{G}_0), \label{con3:prob3} \\
    & \mbox{Eqs.~\eqref{con3:prob1}--\eqref{con6:prob1}}. \label{con4:prob3}
\end{align}
A reasonable choice is to set $\mathcal{G}_0 \coloneqq \{0\}$ with $\bm{o}(0) \coloneqq (1,2,\ldots,R)$, which is the scheduled order of ship arrivals. 
In this case, the corresponding problem size is clearly independent of the number of scenarios.   
Our algorithm repeatedly solves the relaxed problem~\eqref{obj:prob3}--\eqref{con4:prob3} while adding the relaxed constraints (Eqs.~\eqref{con1:prob2}--\eqref{con3:prob2}) to the problem.

We first set $k \leftarrow 1$ and find an optimal solution $(\bar{u}^{(k)}, \bar{\bm{x}}^{(k)}, \bar{\bm{z}}^{(k)}, \bar{\gamma}^{(k)})$ to the relaxed problem~\eqref{obj:prob3}--\eqref{con4:prob3}. 
Its objective value is employed as a lower bound ${\rm LB}$ on problem~\eqref{obj:prob2}--\eqref{con4:prob2}. 
We can also find an upper bound ${\rm UB}$ by substituting the obtained solution $\bar{\bm{x}}^{(k)}$ into problem~\eqref{obj:prob1}--\eqref{con8:prob1}. 

We next develop a new scenario subset as 
\begin{align}
  \mathcal{G}_k \coloneqq \left\{j \in \mathcal{J} ~\middle|~ \sum_{s\in \mathcal{S}}\sum_{h \in \mathcal{H}\setminus\{1\}}\bar{z}^{(k)}_{shj}-\bar{\gamma}^{(k)} > 0 \right\} \label{eq:J_k}.
\end{align}
We then update the relaxed problem~\eqref{obj:prob3}--\eqref{con4:prob3} by adding Eq.~\eqref{con1:prob2} for $\mathcal{G} = \mathcal{G}_k$ and Eqs.~\eqref{con2:prob2}--\eqref{con3:prob2} for $j \in \mathcal{G}_k \setminus (\cup_{\ell=0}^{k-1} \mathcal{G}_{\ell})$ to the problem. 
This process is repeated until the optimality gap is sufficiently small as 
\[
\frac{{\rm UB} - {\rm LB}}{\rm UB} \le \varepsilon, 
\]
where $\varepsilon \ge 0$ is a threshold for optimality. 

Our algorithm is summarized in Algorithm~\ref{alg:alg1}. 
This algorithm is known to terminate within a finite number of iterations~\cite{takano2015cutting,kunzi2006computational,kobayashi2021bilevel}; see K{\"u}nzi-Bay and Mayer~\cite{kunzi2006computational} for details of theoretical properties of the cutting-plane method. 
\begin{algorithm}[H] 
  \caption{Algorithm based on the cutting-plane method for solving problem~\eqref{obj:prob2}--\eqref{con4:prob2}}
  \label{alg:alg1}
  \begin{algorithmic}[1]
  \STATE $k \leftarrow 1$.
  \STATE $\mathcal{G}_0 \coloneqq \{0\}$,~$\bm{o}(0) \coloneqq (1,2,\ldots,R)$.
  \STATE ${\rm LB} \leftarrow 0, ~{\rm UB} \leftarrow \infty$.
  \WHILE {$\displaystyle \frac{{\rm UB}-{\rm LB}}{\rm UB} \ge \varepsilon$}
  \STATE Solve the relaxed problem~\eqref{obj:prob3}--\eqref{con4:prob3}. 
  \STATE Update ${\rm LB}$ and ${\rm UB}$ if possible. \\
  \STATE Compute $\mathcal{G}_k$ as in Eq.~\eqref{eq:J_k}.
  \STATE Add Eq.~\eqref{con1:prob2} for $\mathcal{G} = \mathcal{G}_k$ and Eqs.~\eqref{con2:prob2}--\eqref{con3:prob2} for $j \in \mathcal{G}_k \setminus (\cup_{\ell=0}^{k-1} \mathcal{G}_{\ell})$ to the relaxed problem. 
  \STATE $k \leftarrow k+1$.
  \ENDWHILE
  \end{algorithmic}
\end{algorithm}

\section{NUMERICAL EXPERIMENT}
\label{sec:exper}
In this section, we evaluate the effectiveness of our container pre-marshalling method through numerical experiments.

\subsection{Synthetic datasets}\label{sec:synthetic}
Table \ref{tab:setting} lists the two types of synthetic datasets: $4 \times 4$-Bay and $8 \times 8$-Bay, used in the numerical experiments, where $S$, $H$, $N$, and $R$ are the numbers of stacks (columns), height levels (rows), containers, and priority classes (or ships), respectively. 
The number $|\mathcal{N}_r|$ of containers was set to the same value for all $r \in \mathcal{R}$ in each dataset. 
Note that the total number of containers was set as $N = \left( S - 1 \right) H$, which ensures that containers can be transformed into any layout regardless of the initial layout~\cite{boge2020robust}. 

\begin{table}[h]
    \caption{Synthetic datasets}
    \label{tab:setting}
    \centering
    \scalebox{0.95}{
    \begin{tabular}{cccccc}
    \hline
    Dataset & $S$ & $H$ & $N$ & $R$ & $|\mathcal{N}_r|$ \\
    \hline
    $4 \times 4$-Bay & 4 & 4 & 12 & 6 & 2 \\
    $8 \times 8$-Bay & 8 & 8 & 56 & 14 & 4 \\
    \hline
    \end{tabular}
    }
\end{table}

We sampled ship arrival times from a multivariate normal distribution as 
$\bm{t}(i) \sim \mathrm{N}(\bm{\mu}, \bm{\Sigma})$ for $i \in \mathcal{I}$, where $\bm{\mu} \coloneqq (\mu_r)_{r \in \mathcal{R}}$ and $\bm{\Sigma} \coloneqq (\sigma_{rq})_{(r,q) \in \mathcal {R} \times \mathcal{R}}$ are the mean vector and the covariance matrix, respectively.  
The mean vector was drawn from a uniform distribution on the interval $[0,R]$ as $\mu_r \sim \mathrm{U}(0,R)$ for $r \in \mathcal{R}$. 
The (positive semidefinite) covariance matrix was drawn from a Wishart distribution as $\bm{\Sigma} \sim \mathrm{W}(R, \left(1/R \right)\bm{I})$, where $\bm{I}$ is the identity matrix.
Here, the degrees of freedom and the scale matrix were set to $R$ and $\left(1/R \right)\bm{I}$, respectively. 
In this case, the mean matrix was $\mathrm{E}(\bm{\Sigma}) = \bm{I}$, and the variance of each entry was $\mathrm{Var}(\sigma_{rq}) = 2/R$ if $r=q$, and $\mathrm{Var}(\sigma_{rq}) = 1/R$ otherwise. 

\begin{table*}[t]
  \caption{The number of misplaced containers for the testing scenario set}
  \label{tab:result_place_IW}
  \small
  \centering
   \begin{tabular}{cccrrrrrrrr}
    \hline
     Dataset & \multicolumn{2}{c}{Method}  & \multicolumn{2}{c}{Average} & \multicolumn{2}{c}{75\%-quantile} & \multicolumn{2}{c}{90\%-quantile} & \multicolumn{2}{c}{99\%-quantile} \\
    
    \hline
    $4 \times 4$-Bay & CV@R & $\alpha=0.00$ & 0.203 & $(\pm0.115)$ & $\bm{0.200}$ & $(\pm0.200)$ & 0.800 & $(\pm0.442)$ & $\bm{1.400}$ & $(\pm0.521)$ \\
    & ($|\mathcal{I}|=10^4$) & $\alpha=0.75$ & 0.183 & $(\pm0.116)$ & $\bm{0.200}$ & $(\pm0.200)$ & $\bm{0.600} $ & $(\pm0.427)$ & $\bm{1.400}$ & $(\pm0.521)$\\
    & &$\alpha=0.90$ & $\bm{0.128}$ & $(\pm 0.063)$ & 0.400 & $(\pm0.267)$ & $\bm{0.600}$ & $(\pm0.427)$ & $\bm{1.400}$ & $(\pm0.427)$ \\ \cline{2-11}
    & RO & $\Gamma=0$ & 1.411& $(\pm0.166)$ & 1.900 & $(\pm0.233)$ & 2.800 & $(\pm0.249)$ & 4.800 &  $(\pm0.359)$ \\
    &  & $\Gamma=1$ & 0.725& $(\pm0.192)$ & 1.500 & $(\pm0.453)$ & 2.100 & $(\pm0.526)$ & 3.200 & $(\pm0.814)$ \\
    &  & $\Gamma=2$ & 0.286 & $(\pm0.096)$ & 0.500 & $(\pm0.269)$ & 1.100 & $(\pm0.379)$ & 1.700 & $(\pm0.423)$ \\
    &  & $\Gamma=3$ & 0.231 & $(\pm0.080)$ & $\bm{0.200}$ & $(\pm0.200)$ & 0.800 & $(\pm0.327)$ & $\bm{1.400}$ & $(\pm0.427)$ \\
    &  & $\Gamma=4$ & 0.322 & $(\pm0.080)$ & 0.500 & $(\pm0.167)$ & 0.900 & $(\pm0.180)$ & 1.700 & $(\pm0.367)$ \\ \hline
    $8 \times 8$-Bay & CV@R & $\alpha=0.00 $ & 0.014 & $(\pm0.007)$ & $\bm{0.000}$ & $(\pm0.000)$ & $\bm{0.000}$ & $(\pm0.000)$ & 0.800 & $(\pm0.533)$ \\
    & ($|\mathcal{I}|=10^4$) & $\alpha=0.75$ & $\bm{0.010}$ & $(\pm0.005)$ & $\bm{0.000}$ & $(\pm0.000)$ & $\bm{0.000}$ & $(\pm0.000)$ & $\bm{0.400}$ & $ (\pm0.400)$ \\
    &  & $\alpha=0.90$ & 0.011 & $(\pm0.005)$ & $\bm{0.000} $ & $(\pm0.000)$ & $\bm{0.000}$ & $(\pm0.000)$ & $\bm{0.400}$ & $(\pm0.400)$ \\ \cline{2-11}
    & RO & $\Gamma=0$ & 8.735 & $ (\pm 0.388)$ & 11.300 & $ (\pm 0.448)$ & 13.600  & $ (\pm 0.452)$ & 18.100 & $ (\pm 0.547)$ \\
    & & $\Gamma=2$ & 0.698 & $(\pm0.340)$ & 0.800 & $(\pm0.611)$ & 2.000 & $(\pm1.033)$ & 4.700 & $(\pm1.453)$ \\
    &  & $\Gamma=4$ & 0.141 & $(\pm0.085)$ & $\bm{0.000}$ & $(\pm0.000)$ & 0.600 & $(\pm0.400)$ & 2.000 & $(\pm1.033)$ \\
    &  & $\Gamma=6$ & 0.016 & $(\pm0.013)$ & $\bm{0.000} $ & $(\pm0.000)$ & $\bm{0.000}$ & $(\pm0.000)$ & $\bm{0.400}$ & $(\pm0.400)$ \\
    &  & $\Gamma=8$ & 0.223 & $(\pm0.080)$ & 0.300 & $(\pm0.300)$ & 1.000 & $(\pm0.447)$ & 2.000 & $(\pm0.558)$ \\
    \hline
   \end{tabular}
 \end{table*}

\subsection{Experimental Design}

We compare the performance of the following optimization models for container pre-marshalling under uncertainty of ship arrival times: 
\begin{description}
\item[CV@R:] our optimization model~\eqref{obj:prob1}--\eqref{con8:prob1}; 
\item[RO:] the robust optimization model~\cite{boge2020robust} based on the uncertainty set of ship arrival orders.
\end{description}

The experiments were performed on a computer with an Intel(R) Core(TM) i5-8250U CPU @ 1.60GHz, 8GB RAM, and 4 cores. 
Optimization problems were solved using the mathematical optimization solver Gurobi Optimizer\footnote{\url{https://www.gurobi.com/}} version 10.0.2.

In the CV@R model, the probability level was set as $\alpha \in \{0.00, 0.75, 0.90\}$, where $\alpha = 0.00$ corresponds to the expectation minimization as in prior studies~\cite{zweers2020optimizing,maniezzo2021stochastic}. 
The number of instances of ship arrival times sampled for scenario generation was set as $|\mathcal{I}| \in \{10^2, 10^3, 10^4\}$; see also Section~\ref{sec:sce}. 

In the RO model, the swap distance for specifying the uncertainty set was set as $\Gamma \in \{0,1,2,3,4\}$ for the $4 \times 4$-Bay dataset and $\Gamma \in \{0,2,4,6,8\}$ for the $8 \times 8$-Bay dataset, where $\Gamma = 0$ corresponds to assuming that all ships arrive as scheduled. 

We evaluate the quality of container layouts according to the following process:
\begin{enumerate}
\item Randomly sample ship arrival times as $\bm{t}(i) \sim \mathrm{N}(\bm{\mu}, \bm{\Sigma})$, where the sample size was $|\mathcal{I}|$ for training and $10^4$ for testing; 
\item Compute optimal container layouts using the training scenario set of ship arrival orders; 
\item Calculate the number of misplaced containers from the testing scenario set of ship arrival orders. 
\end{enumerate}

We repeated this process ten times with different values of $(\bm{\mu},\bm{\Sigma})$ and show average values over the ten trials in Tables~\ref{tab:result_place_IW} and \ref{tab:calc_time_cvar}, with standard errors in parentheses. 
We terminated a computation if it did not complete within 3600 s; in those cases, the best feasible solution found within 3600 s was taken as the result.

\subsection{Evaluation of the number of misplaced containers} \label{subsec:evaluation_place}
First, we compare the quality of container layouts determined by the CV@R and RO models.
Table \ref{tab:result_place_IW} shows the average, 75\%-quantile, 90\%-quantile, and 99\%-quantile of the number of misplaced containers for the testing scenario set, where the training sample size was $|\mathcal{I}|=10^4$. 

It is clear from Table \ref{tab:result_place_IW} that our CV@R model outperformed the RO model in the average number of misplaced containers. 
For the $4 \times 4$-Bay dataset, the CV@R model significantly reduced the number of misplaced containers compared to the RO model.
For the $8 \times 8$-Bay dataset, the number of misplaced containers was easily reduced due to the high degree of freedom of container moves, resulting in only a small difference between our CV@R model and the best RO model with $\Gamma = 6$. 

The CV@R model performed better with $\alpha \in \{0.75, 0.90\}$ than with $\alpha = 0$. 
This suggests that CV@R minimization is more effective than expectation minimization in container pre-marshalling under uncertainty. 
Additionally, the performance of the RO model with $\Gamma = 0$ was extremely bad. 
This emphasizes the importance of considering the uncertainty of ship arrival times in the container pre-marshalling problem.

The CV@R model treats uncertain ship arrival times as a multivariate probability distribution, whereas the RO model is unable to exploit such a multivariate relationship of ships arrival times. 
For this reason, the CV@R model is superior to the RO model in the quality of container layouts under uncertainty. 

\begin{table*}[t]
\caption{Computation time required for solving optimization models}
\small
\centering
\label{tab:calc_time_cvar}
\begin{tabular}{cccrrccccrr}
\hline 
\multicolumn{5}{c}{$4 \times 4$-Bay dataset} & & \multicolumn{5}{c}{$8 \times 8$-Bay dataset } \\
\cline{1-5} \cline{7-11} 
\multicolumn{2}{c}{ Method } & Solver & \multicolumn{2}{c}{Time} & &  \multicolumn{2}{c}{ Method } & Solver & \multicolumn{2}{c}{Time} \\
\hline
CV@R & $|\mathcal{I}|=10^2$ & Algorithm~\ref{alg:alg1} & 1.50 & $(\pm0.17)$ & & CV@R & $|\mathcal{I}|=10^2$ & Algorithm~\ref{alg:alg1} & 22.80 & $(\pm4.73)$ \\
($\alpha=0.75$) & & Gurobi & 1.84 & $(\pm0.80)$ & &  ($\alpha=0.75$) & & Gurobi & $>3600$ & \\
\cline{2-5}\cline{8-11}
 & $|\mathcal{I}|=10^3$ & Algorithm~\ref{alg:alg1} & 8.29 & $(\pm0.75)$ & &  & $|\mathcal{I}|=10^3$ & Algorithm~\ref{alg:alg1} & 723.46 & $(\pm35.72)$ \\
 & & Gurobi & 6.57 & $(\pm2.88)$ & & & & Gurobi & \multicolumn{2}{c}{Memory error} \\
\cline{2-5}\cline{8-11}
 & $|\mathcal{I}|=10^4$ & Algorithm~\ref{alg:alg1} & 38.47 & $(\pm1.77)$ & & & $|\mathcal{I}|=10^4$ & Algorithm~\ref{alg:alg1} & $>3600$ & \\
 & & Gurobi & 12.99 & $(\pm4.63)$ & & & & Gurobi & \multicolumn{2}{c}{Memory error} \\
\cline{2-4}\cline{7-9}
\hline
RO & $\Gamma=0$ & Gurobi & 0.05 & $(\pm0.00)$ & & RO & $\Gamma=0$ & Gurobi & 0.08 & $(\pm0.00)$\\
 & $\Gamma=1$ & Gurobi & 0.03 & $(\pm0.00)$ & & & $\Gamma=2$ & Gurobi & 0.10 &$(\pm0.01)$ \\
 & $\Gamma=2$ & Gurobi & 0.04 & $(\pm0.00)$ & & & $\Gamma=4$ & Gurobi & 0.21 & $(\pm0.01)$ \\
 & $\Gamma=3$ & Gurobi & 0.04 & $(\pm0.00)$ & & & $\Gamma=6$ & Gurobi & 0.42 & $(\pm0.02)$ \\
 & $\Gamma=4$ & Gurobi & 0.04 & $(\pm0.00)$ & & & $\Gamma=8$ & Gurobi & 117.10 & $(\pm2.25)$ \\
 \hline
\end{tabular}
\end{table*}

\subsection{Evaluation of the computation time}
Next, we evaluate the computational efficiency of our algorithm (Algorithm~\ref{alg:alg1}), compared with direct application of the optimization solver (Gurobi) to problem~\eqref{obj:prob1}--\eqref{con8:prob1}. 
Table \ref{tab:calc_time_cvar} shows the computation time required for solving optimization models, where the probability level was $\alpha = 0.75$. 
Note that ``$>3600$'' means that the computation time exceeded 3600 s in any of the ten trials. 

For the CV@R model, the computation times of both Algorithm~\ref{alg:alg1} and Gurobi markedly increased as the sample size ($|\mathcal{I}|$) increased. 
For the $4 \times 4$-Bay dataset, Gurobi was always faster than Algorithm~\ref{alg:alg1}. 
For the $8 \times 8$-Bay dataset, however, Gurobi terminated the computation due to the time limit when $|\mathcal{I}| = 10^2$ and memory error when $|\mathcal{I}| \in \{10^3, 10^4\}$. 
On the other hand, Algorithm~\ref{alg:alg1} finished the computation within the time limit when $|\mathcal{I}| \in \{10^2, 10^3\}$ and avoided the memory error even when $|\mathcal{I}| = 10^4$. 
These results suggest that our algorithm is very effective in computing solutions to large-scale container pre-marshalling problems under uncertainty.

The computation time of the RO model was almost independent of $\Gamma$ for the $4 \times 4$-Bay dataset, whereas the computation time largely increased with $\Gamma$ for the $8 \times 8$-Bay dataset. 
For both datasets, the RO model was faster than the CV@R model, even though the RO model was inferior to the CV@R model in terms of the container layout quality. 

\section{CONCLUSION}
\label{sec:concl}
We considered a container pre-marshalling problem under uncertainty of ship arrival times. 
For this problem, we proposed a mixed-integer linear optimization model for minimizing CV@R for the number of misplaced containers, where the uncertainty of ship arrival times was represented as multiple scenarios generated from a multivariate probability distribution.
Moreover, we developed an exact algorithm based on the cutting-plane method to handle large-scale problems.

Through numerical experiments using synthetic datasets, we confirmed that our CV@R model is more effective than the robust optimization model~\cite{boge2020robust} in terms of the quality of container layouts under uncertainty. 
In addition, our exact algorithm was highly effective in computing solutions to large-scale container pre-marshalling problems.

A future direction of study will be to develop an algorithm for solving the container pre-marshalling problem faster. 
To this end, we can use speed-up techniques~\cite{takano2015cutting,kobayashi2021bilevel} of the cutting-plane method or design heuristic algorithms specialized for container pre-marshalling problems. 
Another direction of future research will be to simultaneously determine a sequence of container moves and a final container layout under uncertainty. 

\bibliographystyle{cas-model2-names}

\bibliography{cite}


\end{document}